\documentclass[12pt]{article}
\def\q{\hfill\rule{1ex}{1ex}}
\def\0{\emptyset}

\def\p{{\bf Proof.} \quad}
\def\q{\hfill\rule{1ex}{1ex}}

\def\n{\noindent}

\usepackage{graphicx}
 \usepackage{amsmath}
 \usepackage{amsfonts}

 \usepackage{algorithm}
 \usepackage{algorithmic}

\begin{document}

\title{\bf $w$-Dominating Set Problem on Graphs of Bounded Treewidth }
\author{{\small\bf Ke Liu}\thanks{email:  liuke17@mails.tsinghua.edu.cn
}\quad {\small\bf Mei
Lu}\thanks{email: lumei@tsinghua.edu.cn}\\
{\small Department of Mathematical Sciences, Tsinghua
University, Beijing 100084, China.}}

\date{}

\maketitle\baselineskip 16.3pt

\begin{abstract}
Let $G=(V,E)$ be a graph.    Let $w$ be a positive integer. A $w$-dominating set is a vertex subset $S$ such that for all  $v\in V$, either $v\in S$ or it has at least $w$ neighbors in $S$. The $w$-Dominating Set problem is to find the minimum $w$-dominating set. The $L$-Max $w$-Dominating Set problem is to find the vertex subset $S$ of cardinality at most $L$ that maximizes $|S|+|\{v\in V\setminus S~|~|N(v)\cap S|\geq w\}|$, where $N(v)=\{u|uv\in E\}$. In this paper, we   give  polynomial time algorithms to $w$-Dominating Set problem and $L$-Max $w$-Dominating Set problem
 on graphs of bounded treewidth.
\end {abstract}


\vskip.3cm

\n{\large\bf 1.\ Introduction}

Let $G=(V,E)$ be a graph. A {\em dominating set} is a subset $S$ of $V$ such that for all $v\in V$, either $v\in S$ or it has a neighbor in $S$. The size of the minimum dominating is called the {\em domination number}. Let $w$ be a positive integer. A $w$-dominating set is a vertex subset $S$ such that for all  $v\in V$, either $v\in S$ or it has at least $w$ neighbors in $S$. The $w$-Dominating Set problem is to find the minimum $w$-dominating set. The $L$-Max $w$-Dominating Set problem is to find the vertex subset $S$ of cardinality at most $L$ that maximizes $|S|+|\{v\in V\setminus S~|~|N(v)\cap S|\geq w\}|$, where $N(v)=\{u|uv\in E\}$.

Dominating set problem is shown to be W[2]-complete \cite{Nie}, but it is solvable in polynomial time on graphs of bounded treewidth. Let $tw$ be the treewidth of a graph. Alber et al. \cite{AlNi} gave a $4^{tw}n^{O(1)}$ time algorithm for the dominating set problem, improving  the  $9^{tw}n^{O(1)}$ algorithm given by Telle et al. \cite{Telle}. Using fast subset convolution (see Bjorklund et al. \cite{BCo} or Cygan et al. \cite{Cygan}), the running time of the algorithm in \cite{AlNi} can be improved to $3^{tw}n^{O(1)}$. Roayaei et al. gave a $4^{tw}L^2n^{O(1)}$ time algorithm \cite{RImpro} for the  $L$-Max Dominating Set problem.

Dinh et al. studied the Positive Influence Dominating Set (PIDS) in \cite{Dinh}. Let $G=(V,E)$ be a graph. A PIDS is a vertex set such that for all  $v\in V$, either $v$ is selected into PIDS, or it has at least $\rho d(v)$ neighbors in PIDS for some constant $0<\rho<1$. The PIDS problem is to find a PIDS of minimum cardinality. The authors proved the inapproximability factor $(1/2-\epsilon \ln n)$ for PIDS problem and proposed a linear-time algorithm to find optimal solutions of PIDS over trees.
Note that if $G$ is regular, then the Positive Influence Dominating Set problem is just the $w$-Dominating Set problem.

The treewidth of a graph is an important invariant in structural and algorithmic graph theory. The concept of treewidth was originally introduced by Bertel\'e et al. \cite{BB} under the name of dimension. It was later rediscovered by Halin \cite{Hal} in 1976 and by  Robertson  et al. \cite{RS84} in 1984, respectively. Now it has  been studied by many other authors (see for example \cite{GM}-\cite{KS}).
 The treewidth of a graph gives an indication of how far away the graph is from being
a tree or forest. The closer the graph is to being a forest, the smaller is its treewidth.
 The treewidth
is a graph parameter that plays a fundamental role in various graph algorithms. It is well-known that many NP-complete problems can be solvable in polynomial time on graphs of bounded treewidth \cite{Bo}. In this paper, we will consider the following two problems.

    {\bf $w$-Dominating-Set}

    {\em Instance:} Graph $G=(V,E)$, $k,L\in \mathbb{N}$.

    {\em Parameter:} treewidth $k$.

    {\em Problem:} Decide whether $G$ has a $w$-dominating set of cardinality at most $L$.


    {\bf $L$-Max $w$-Dominating-Set}

    {\em Instance:} Graph $G=(V,E)$, $k,L\in \mathbb{N}$.

    {\em Parameter:} treewidth $k$.

    {\em Problem:} Find a vertex subset $S$ of cardinality at most $L$ such that  $|S|+|\{v\in V\setminus S~|~|N(v)\cap S|\geq w\}|$ as larger as possible.
\vskip.2cm


%
%
%
%
%
%
%
%




Our main results in this paper are to give  polynomial time algorithms to $w$-Dominating Set problem and $L$-Max $w$-Dominating Set problem on the graphs of bounded treewidth, which are shown in Sections 3 and 4, respectively. By setting $w=1$ to the two main results, we get FPT algorithm to Dominating Set problem which is also proved in \cite{Cygan} and a better FPT algorithm to $L$-Max Dominating Set problem than that in \cite{RImpro}.
\vskip.2cm

\n{\large\bf 2.\ Treewidth}
\vskip.2cm

In this Section, we give  definitions involving in treewidth.  The treewidth of a graph is defined through
the concept of tree decompositions. In the following, we will use $T$ to denote the vertex set of $T$ when $T$ is a tree and we call the vertex of $T$ {\em node}.

{\bf Definition 2.1 }
A tree decomposition of a graph $G=(V,E)$ is a pair $(T,(B_{t})_{t\in T})$, where $T$ is a tree and $(B_{t})_{t\in T}$ a family of subsets of $V$ such that:

\n (1) for every $v\in V$, the set $B^{-1}(v)=\{t\in T|v\in B_{t}\}$ is nonempty and connected in $T$;

\n (2) for every edge $\{u,w\}\in E$, there is a $t\in T$ such that $u,w\in B_{t}$.

\n The width of the decomposition $(T,(B_{t})_{t\in T})$ is the number
$$\max\{|B_{t}||t\in T\}-1.$$
The treewidth $tw(G)$ of $G$ is the minimum of the widths of the tree decompositions of $G$.
\vskip.2cm
By Definition 2.1, each graph $G=(V,E)$ has a tree decomposition $(T,(B_{t})_{t\in T})$ where $T$ contains only one node $t$ with $B_{t}=V$. And this kind of decomposition has width $|V|-1$ which is the largest width of the graphs on $|V|$ vertices. A {\em rooted tree decomposition} is a tree decomposition with a distinguished root
node, denoted by $r$. Given a rooted tree decomposition $(T,(B_{t})_{t\in T})$ with a root node $r$ and a node $t$ of $T$, let $Desc(t)$
be the set of descendants of node $t$ in $T$, including $t$; let $T_t = T [Desc(t)]$ be a subtree of $T$ rooted at $t$; let $G_t =
G[\cup_{s\in T_t}B_s]$. Then $T_r=T$ and $G_r=G$.

In order to do algorithm analysis and reduce the time complexity on tree decomposition, we need the definition of nice tree decomposition \cite{ABo}.
\vskip.2cm

{\bf Definition 2.2 }A tree decomposition $(T,(B_{t})_{t\in T})$ of a graph $G$ is {\em nice} if it satisfies the following properties:

\n (1) every node of $T$ has at most two child nodes;

\n (2) if a node $t\in T$ has two child nodes $t_1,t_2$, then $B_t=B_{t_1}=B_{t_2}$ and $t$ is called a {\em join node};

\n (3) if a node $t\in T$ has one child node $t_1$, then one of the following must hold

(a) $B_t=B_{t_1}\cup \{x_0\}$ for a $x_0\in V(G)$, $t$ is called an {\em introduce node}, or

(b) $B_t=B_{t_1}\setminus \{x_0\}$ for a $x_0\in V(G)$, $t$ is called a {\em forget node}.
 \vskip.2cm

The number of nodes of a nice tree decomposition can be controlled by $n^{O(1)}$, where $n=|V(G)|$. And it is not hard to transform a given tree decomposition into a nice one \cite{Knice}. The problem of deciding whether a graph has tree decomposition of treewidth at most $k$ is NP-complete \cite{A} and Bodlaender \cite{Bo} proved that the problem is fixed-parameter tractable.

{\bf Lemma 2.1 \cite{Bo}}
There is a polynomial $p$ and an algorithm that, given a graph $G=(V,E)$, computes a tree decomposition of $G$ of width $tw(G)=k$ in time at most
$2^{p(k)}|V|$.

\vskip.2cm

By Lemma 2.1, we can get a nice tree decomposition of $G$ in polynomial time when $tw(G)$ equals to a constant $k$. Thus in the rest of our paper, we assume we have a nice tree decomposition of $G$ with width $tw(G)$.
\vskip.2cm
We also need the definition of subset convolution. Given a set $S$ and two functions $g,h: 2^S\rightarrow \mathbb{Z}$, the subset convolution of $g$ and $h$ is a function
$(g*h):2^S\rightarrow \mathbb{Z}$ such that for $Y\subseteq S$,
$$(g*h)(Y)=\min\limits_{\substack{A\cup B=Y\\ A\cap B=\emptyset}}(g(A)+h(B))$$
or
$$(g*h)(Y)=\max\limits_{\substack{A\cup B=Y\\ A\cap B=\emptyset}}(g(A)+h(B)).$$

For the complexity of computing subset convolution, we have the following result.

{\bf Lemma 2.2 \cite{BCo}}
Let $S$ be a set with $n$ elements and $M$ be a positive integer. For two functions $g,h:2^S\rightarrow\{-M,\ldots,M\}\cup\{+\infty\}$, if all the values of $g$ and $h$ are given, then all the $2^n$ values of the subset convolution of $g$ and $h$ can be computed in $2^nn^{O(1)}O(M\log(Mn)\log\log(Mn))$ time.

\vskip.2cm
\n{\large\bf 3.\ $w$-Dominating Set Problem}
\vskip.2cm
In this Section, we study $w$-Dominating Set problem on graphs with bounded treewidth.

\vskip.2cm

{\bf Theorem 3.1 }
{\em Let $k,n,w$ be three positive integers and $G=(V,E)$ be a graph of order $n$  with $tw(G)=k$. Then the $w$-Dominating Set problem can be solved in time
$(\frac{(w+1)(w+2)}{2})^kk^{O(1)}n.$}

\p Let $(T,(B_t)_{t\in T})$ be a nice tree decomposition of $G$ rooted at $r$ with width $k$. Then $|B_t|\le k+1$ for all $t\in T$.
For each bag $B_t$, the coloring of $B_t=\{x_1,\ldots,x_{|B_t|}\}$ is a mapping $f_t:B_t\rightarrow \{0,1,\ldots,w,+\infty\}$ assigning $w+2$ different colors to the vertices in the bag, and the color assigned to the vertex $x$ is denoted by $f_t(x)$. We use  a vector $(f_t(x_1),\ldots,f_t(x_{|B_t|}))$ to denote a coloring of $B_t$, that is, $f_t=(f_t(x_1),\ldots,f_t(x_{|B_t|}))$. There exists at most $(w+2)^{k+1}$ colorings of $B_t$.

For a coloring $f_t$ of $B_t$, let $Dom(t,f_t)\subseteq \cup_{t'\in T_t}B_{t'}$ denote the minimum vertex set subject  to (i) $|N(x)\cap Dom(t,f_t)|\geq w$ for all $x\in (\cup_{t'\in T_t}B_{t'})\setminus (B_t\cup Dom(t,f_t))$; (ii) $Dom(t,f_t)\cap B_t=f_t^{-1}(+\infty)$; (iii) $|N(x)\cap Dom(t,f_t)|\geq f_t(x)$ for all $x\in B_t\setminus f^{-1}_t(+\infty)$. The evaluation index $c[t,f_t]$ is defined as
$$c[t,f_t]=\left\{
\begin{array}{ll}
+\infty  &   \mbox{there is no such set } Dom(t,f_t),\\
 |Dom(t,f_t)|  &  \mbox{otherwise.}
\end{array}\right.$$
Then when we get all the $c[r,f_r]$, we have solved the $w$-Dominating Set problem.

To make the process of calculation of $c[t,f_t]$ clear, we need the following definition of the partial ordering of coloring $f_t$. First, on the color set $\{0,1,\ldots,w,+\infty\}$, let $\prec$ denote the partial ordering defined by two basic rules:

\n (1) $i\prec j$ for $i,j\in\{0,1,\ldots,w\}$ if $i\leq j$;

\n (2) $i\prec i$ for $i\in\{0,1,\ldots,w,+\infty\}$.

Then, we say $f_t\prec f'_t$ if and only if $f_t(x)\prec f'_t(x)$ for all $x\in B_t$. Based on the partial ordering, we know $f_t\prec f_t'$ implies that $c[t,f_t]\leq c[t,f'_t]$.

When $t$ is a leaf node,  we have $T_t=t$. For each coloring $f_t$, we know
$$c[t,f_t]=\left\{
\begin{array}{ll}
+\infty  &   \exists x\in B_t\setminus f_t^{-1}(+\infty),\left|N(x)\cap f_t^{-1}(+\infty)\right|<f_t(x),\\
 \left|f_{t}^{-1}(+\infty)\right|  &  \mbox{otherwise.}
\end{array}\right.$$

It takes $O((w+2)^{k+1}(k+1))$ time to compute all $c[t,f_t]$ for each leaf node $t$. After calculating the evaluation index of $f_t$ for leaf nodes, we visit the bags of the tree decomposition from leaves to the root and calculate the corresponding evaluation index in each step according to the following rules.

Let $t$ be a non-leaf node and assume we have count all $c[t',f_{t'}]$, where $t'$ is a child node of $t$. We then design algorithm for three types of $t$.

\textbf{Forget node:} Suppose $t$ is a forget node. Assume $B_{t'}=\{x_1,\ldots,x_{|B_{t'}|},x_0\}$ and $B_t=B_{t'}\setminus \{x_0\}$ for a $x_0\in V(G)$. For each $f_t$, define $f_{t'}=f_t\times \{d\}=(f_t(x_1),\ldots,f_t(x_{|B_t|}),d)$, where $f_{t'}(x_0)=d$ and $d\in \{0,1,\ldots,w,+\infty\}$. If $d\not=+\infty$, then $x_0\notin Dom(t',f_{t'})$. Since $B^{-1}(x_0)$ is connected in $T$ and $x_0\notin B_t$, $N(x_0)\subseteq \cup_{u\in T_{t'}}B_{u}$  which implies $d=w$. Thus
we have
$$c[t,f_t]=\min\limits_{d\in \{+\infty ,w\}}c[t',f_t\times \{d\}].$$
It takes $O((w+2)^{k+1}(k+1))$ time for each forget node.

\textbf{Introduce node:} Suppose $t$ is an introduce node. Assume $B_{t'}=\{x_1,\ldots,x_{|B_t'|}\}$ and $B_t=B_{t'}\cup \{x_0\}$ for a $x_0\in V(G)$. For each  $f_t=(f_{t}(x_1),\ldots,f_{t}(x_{|B_{t'}|}),f_t(x_0))$, we have $f_t=f_{t'}\times \{f_t(x_0)\}$, where $f_{t'}=(f_{t}(x_1),\ldots,f_{t}(x_{|B_{t'}|}))$.
 Define $f'_{t'}=(f'_{t'}(x_1),\ldots,f'_{t'}(x_{|B_{t'}|}))$ such that for all $y\in B_{t'}$,
$$f'_{t'}(y)=\left\{
\begin{array}{ll}
+\infty  &   \mbox{if } f_{t}(y)=+\infty, \\
\max(0,f_{t}(y)-1)  &  \mbox{if } f_{t}(y)\neq +\infty \mbox{ and } y\in N(x_0),\\
f_{t}(y)           &  \mbox{otherwise.}
\end{array}\right.$$
Then $f'_{t'}$ is a coloring of $B_{t'}$. We claim the following results hold:

\n (3) $c[t,f_{t'}\times \{+\infty\}]=c[t',f'_{t'}]+1$;

\n (4) $c[t,f_{t'}\times \{s\}]=c[t',f_{t'}]$ if $|N(x_0)\cap f^{-1}_{t'}(+\infty)|\geq s$, where $0\leq s\leq w$;

\n (5) $c[t,f_{t'}\times \{s\}]=+\infty$ if $|N(x_0)\cap f^{-1}_{t'}(+\infty)|< s$, where $0\leq s\leq w$.

If  $f_t(x_0)=+\infty$,  $x_0$ is selected into $Dom(t,f_t)$. Then  $Dom(t',f'_{t'})\cup\{x_0\}=Dom(t,f_t)$ and we have (3).

If   $f_t(x_0)=s$ for  $0\leq s\leq w$, we need to check whether  $|N(x_0)\cap f^{-1}_{t'}(+\infty)|\ge s$. Since $B^{-1}(x_0)$ is connected in $T$, $N(x_0)\cap ((\cup_{u\in T_{t'}}B_{u})\setminus B_{t'})=\0$.  Thus (4) and (5) hold.

Since we need $O(k+1)$ time to calculate the number of neighbors selected into the $w$-dominating set, it takes $O((w+2)^{k+1}(k+1))$ time for each introduce node.

\textbf{Join node:} Suppose $t$ is a join node. Assume its child nodes are $t_1,t_2$ with $B_t=B_{t_1}=B_{t_2}$. For each   $f_t$, let $s'(u)=|N(u)\cap f_t^{-1}(+\infty)|$ for $u\in B_t\setminus f_t^{-1}(+\infty)$. We call $\{f_{t_1},f_{t_2}\}$ {\em a good pair} of $f_t$ if

\n (6) $f_t(x)=+\infty$ if and only if $f_{t_1}(x)=f_{t_2}(x)=+\infty$ for any $x\in B_t$;

\n (7) for any $0\leq s\leq w$ and $x\in B_t$, $f_t(x)=s$ if and only if $f_{t_1}(x)+f_{t_2}(x)-s'(x)\geq s.$

Obviously, for each $f_t$, such good pair exists.
Given $f_t$, let $\{f_{t_1},f_{t_2}\}$ be a good pair of $f_t$. Then
$$Dom(t,f_t)=Dom(t_1,f_{t_1})\cup Dom(t_2,f_{t_2}).$$
Since $(\cup_{t'\in T_{t_1}}B_{t'})\cap (\cup_{t'\in T_{t_2}}B_{t'})\subseteq B_t$, $Dom(t_1,f_{t_1})\cap Dom(t_2,f_{t_2})=f_{t}^{-1}(+\infty)$. So  we have
$$|Dom(t,f_t)|=|Dom(t_1,f_{t_1})|+|Dom(t_2,f_{t_2})|-|f_{t}^{-1}(+\infty)|$$
and
$$c[t,f_t]=\min\limits_{\{f_{t_1},f_{t_2}\} \mbox{~is a good pair}} \left\{c[t_1,f_{t_1}]+c[t_2,f_{t_2}]-|f_{t}^{-1}(+\infty)|\right\}.$$
Let $\{f_{t_1},f_{t_2}\}$ and $\{f'_{t_1},f'_{t_2}\}$ be good pairs of $f_t$. For any $x\in B_t$ with $0\leq f_t(x)\leq w$, suppose $f_{t_1}(x)+f_{t_2}(x)-s'(x)= f_t(x)$ and $f_{t_1}(x),f_{t_2}(x)\geq s'(x)$. Then we have $f_{t_1}\prec f'_{t_1}$ and $f_{t_2}\prec f'_{t_2}$ which implies $c[t_1,f_{t_1}]+c[t_2,f_{t_2}]-|f_{t}^{-1}(+\infty)|\le c[t_1,f'_{t_1}]+c[t_2,f'_{t_2}]-|f_{t}^{-1}(+\infty)|$. Hence we can replace (7) by



\n ($7'$)  for any $0\leq s\leq w$ and $x\in B_t$, $f_t(x)=s$  if and only if $f_{t_1}(x)+f_{t_2}(x)-s'(x)= s$ and $f_{t_1}(x),f_{t_2}(x)\geq s'(x)$.

By ($7'$), if $f_t(x)=s$ for some $0\leq s\leq w$, then $(f_{t_1}(x),f_{t_2}(x))\in \{(s'(x),s),(s'(x)+1,s-1)\ldots,(s'(x)+i,s-i),\ldots,(s,s'(x))\}$. Thus (6) and ($7'$)  are equivalent to the following three conditions:

\n (a) $f_t^{-1}(+\infty)=f_{t_1}^{-1}(+\infty)=f_{t_2}^{-1}(+\infty)$;


\n (b) $f_t^{-1}(s)=\{x\in B_t~|~f_{t_1}(x)=s_1,f_{t_2}(x)=s_2,s_1+s_2-s'(x)=s\}$, where $0\leq s\leq w$ and  $s'(x)\leq s_1,s_2\leq s$; particularly,


\n (c) $f_t^{-1}(1)=(f_{t_1}^{-1}(1)\setminus U )\cup (f_{t_2}^{-1}(1) \setminus
U)$ and $(f_{t_1}^{-1}(1)\setminus U)\cap (f_{t_2}^{-1}(1) \setminus
U)=\emptyset,$ where $U=f^{-1}_t (+\infty)\cup f^{-1}_t(2)\cup\cdots\cup f_t^{-1}(w)$.

For disjoint sets $R_{+\infty},R_2,\ldots,R_{w}\subseteq B_t$, let $\overline{R}=R_{+\infty}\cup R_2\cup\ldots\cup R_w$ and $$\mathcal{F}_{\overline{R}}=\{f_t|f_t^{-1}(i)=R_i,i\in \{2,\ldots,w, +\infty\}\}.$$ We want to compute $c[t,f_t]$ for all $f_t\in \mathcal{F}_{R_{+\infty},R_{2},\ldots,R_{w}}$.

For each $f_t\in \mathcal{F}_{\overline{R}}$, if $f_t^{-1}(1)$ is determined, then $f_t^{-1}(0)$ is determined since
 $f_t^{-1}(0)=B_t\setminus (f_t^{-1}(1)\cup \overline{R})$. In this sense, we can think about that $f_t\in \mathcal{F}_{\overline{R}}$ is determined by $f_t^{-1}(1)$. Let $R_1\subset B_t\setminus \overline{R}$, we will rewrite $f_t\in \mathcal{F}_{\overline{R}}$ with $f_t^{-1}(1)=R_1$ by $f_t^{R_1}$.

 Given $R_1\subset B_t\setminus \overline{R}$ and $f_t^{R_1}\in \mathcal{F}_{\overline{R}}$. Let $f_t^{R_1}(x)=s_x$ for all $x\in B_t\setminus R_{+\infty}$, where $2\le s_x\le w$. Denote
$$\mathcal{B}_{R_{1}}=\{\{f_{t_1},f_{t_2}\}|f^{-1}_{t_1}(+\infty)=f^{-1}_{t_2}(+\infty)=R_{+\infty},f_{t_1}(x)+f_{t_2}(x)-s'(x)=s_x \}.$$
Then $|\mathcal{B}_{R_{1}}|\le \prod\limits_{s=2}^{s=w}(s+1)^{|R_s|}$. For any $x\in B_t\setminus \overline{R}$, we have $f_t^{R_1}(x)=0$ or 1. In order to make a pair $\{f_{t_1},f_{t_2}\}\in \mathcal{B}_{R_{1}}$ to be a good pair of $f_t^{R_1}(x)$, we have $f_{t_i}(x)=0$ or 1 for $i\in \{1,2\}$. Thus $f_{t_i}$ is determined by $f_{t_i}^{-1}(1)\cap \overline{R}$ for $i\in \{1,2\}$ and we also rewrite $f_{t_i}$ by $f_{t_i}^{R_1^i}$ if $f_{t_i}^{-1}(1)\cap \overline{R}=R_1^i$. By (c), $f_t^{R_1}(x)=1$ implies that $f_{t_1}(x)=0$ if and only if $f_{t_2}(x)=1$.  Let $R^1_1,R^2_1\subseteq B_t\setminus \overline{R}$ such that $R^1_1\cup R^2_1=R_1$ and $R^1_1\cap R^2_1=\emptyset$. Then the evaluation index of $f_t^{R_1}$ can be calculated by




$$c\left[t,f_t^{R_1}\right]=\min\limits_{\{f_{t_1}^{R_1},f_{t_2}^{R_1}\}\in \mathcal{B}_{R_{1}}}c\left[t,f_t^{R_1},f_{t_1}^{R_1},f_{t_2}^{R_1}\right],$$
where
$$c\left[t,f_t^{R_1},f_{t_1}^{R_1},f_{t_2}^{R_1}\right]=\min\limits_{\substack{R^1_1\cup R^2_1=R_1\\R^1_1\cap R^2_1=\emptyset}} \left(c\left[t_1,f_{t_1}^{R^1_1}\right]+c\left[t_2,f_{t_2}^{ R^2_1}\right]\right)-|R_{+\infty}|.$$
By Lemma 2.2, we can compute $c\left[t,f_t^{R_1},f_{t_1}^{R_1},f_{t_2}^{R_1}\right]$ for every $f_t^{R_1}\in \mathcal{F}_{\overline{R}}$ and a good pair in $2^{|B_t\setminus \overline{R}|}(k+1)^{O(1)}$ time. Then the total complexity of computing  $f_t$ is
$$\begin{array}{rcl}
& &\sum\limits_{R_{+\infty}\subseteq B_t}\left(\binom{|B_t|}{|R_{+\infty}|}\sum\limits_{R_w\subseteq B_t\setminus R_{+\infty}}\left(\binom{|B_t|-|R_{+\infty}|}{|R_w|}(w+1)^{|R_w|}\cdots \right.\right.\\
& &\left.\left.\sum\limits_{R_2\subseteq B_t\setminus (R_{+\infty}\cup R_3\cup\cdots\cup R_w)}\left(\binom{|B_t|-|R_{+\infty}|- \sum_{i=3}^w|R_i|}{|R_2|}3^{|R_2|}2^{|B_t|-|R_{+\infty}|- \sum_{i=2}^w|R_i|}\right)\right)\right)(k+1)^{O(1)}\\
&=&\sum\limits_{R_{+\infty}\subseteq B_t}\left(\binom{|B_t|}{|R_{+\infty}|}(2+3+\ldots+w+1)^{|B_t|-|R_{+\infty}|}\right)(k+1)^{O(1)}\\
&=&(\frac{(w+1)(w+2)}{2})^{|B_t|}(k+1)^{O(1)}.
\end{array}$$

So it takes $(\frac{(w+1)(w+2)}{2})^{k+1}(k+1)^{O(1)}$  time  for each join node. To sum up,  the total time for the algorithm is $(\frac{(w+1)(w+2)}{2})^{k}k^{O(1)}n$.
\q

If we set $w=1$, we get the following corollary.

{\bf Corollary 3.1 \cite{Cygan}}
{\em Let $k,n$ be two positive integers and $G=(V,E)$ be a graph of order $n$ with $tw(G)=k$. Then the Dominating Set problem can be solved in time
$O(3^kk^{O(1)}n).$}


\vskip.2cm

\n{\large\bf 4.\ $L$-Max $w$-Dominating Set problem}
\vskip.2cm

In this section, we consider the $L$-Max $w$-Dominating Set problem on graphs with bounded treewidth.

We are going to give our algorithm  on a nice tree decomposition $(T,(B_{t})_{t\in T})$ with $|T|=n^{O(1)}$.
\vskip.2cm

{\bf Theorem 4.1 }
{\em Let $L,k,n,w$ be four positive integers and $G=(V,E)$ be a graph of order $n$ with $tw(G)=k$. Then the $L$-Max $w$-Dominating Set problem can be solved in  time $(\frac{L(L+1)}{2})(\frac{(w+1)(w+2)}{2})^{k+1}(k+1)^{O(1)}n$.}

\p Let $(T,(B_t)_{t\in T})$ be a nice  tree decomposition of $G$ rooted at $r$ with width $k$. We define $dom(V_2\rightarrow V_1)=\{v\in V_1\setminus V_2~|~|N(v)\cap V_2|\geq w\}$.

For each bag $B_t$, we define the coloring mapping of $B_t$ as Section 3. 
For $f_t$ of $B_t$ and $z\in \{0,1,\ldots,L\}$, we define
$$Dom(t,f_t,z)=arg\max\limits_{S\subseteq \cup_{t'\in T_t}B_{t'}}|dom(S\rightarrow\cup_{t'\in T_t}B_{t'})|+|S|,$$
where $S$ is under the restriction (i) $|S|\leq z$; (ii) $S\cap B_t=f_t^{-1}(+\infty)$; (iii) $|N(x)\cap S|\geq f_t(x)$ for all $x\in B_t\setminus f^{-1}_t(+\infty)$. The evaluation index $c[t,f_t,z]$ is defined as 
$c[t,f_t,z]=-\infty$ if there is no such set $Dom(t,f_t,z)$; otherwise $c[t,f_t,z]=|dom(Dom(t,f_t,z)\rightarrow\cup_{t'\in T_t}B_{t'})|+ |Dom(t,f_t,z))|$.
When we get all the $c[r,f_r,L]$, we  solve the $L$-Max $w$-Dominating Set problem.



For $t\in T$, we count $c[t,f_t,z]$ for a coloring mapping $f_t$ of $B_t$ and $z\in \{0,1,\ldots,L\}$ by induction on the node of $T$.

When $t$ is a leaf node,  we have $T_t=t$. For each  $f_t$, we have $c[t,f_t,z]=-\infty$ if there is $x\in B_t\setminus f_t^{-1}(+\infty)$ such that $|N(x)\cap f_t^{-1}(+\infty)|<f_t(x)$ or $|f_t^{-1}(+\infty)|>z$; otherwise $c[t,f_t,z]=|f_t^{-1}(w)|+|f_t^{-1}(+\infty)|$.

It takes $O((w+2)^{k+1}(k^2+L))$ time for each leaf node. After calculating all $f_t$ for leaf nodes, we visit the bags of the tree decomposition from leaves to the root and calculate the corresponding evaluation index in each step according to the following rules.

Let $t$ be a non-leaf node and assume we have count all $c[t',f_{t'},z]$ for $z\in\{0,1,2,\ldots,L\}$, where $t'$ is a child node of $t$. We  design algorithm for three types of $t$.

\textbf{Forget node:} Suppose $t$ is a forget node. Assume $B_{t'}=\{x_1,\ldots,x_{|B_{t}|},x_0\}$ and $B_t=B_{t'}\setminus \{x_0\}$ for a $x_0\in V(G)$. For each $f_t$ and $z$, we have
$$c[t,f_t,z]=\max\limits_{d\in \{+\infty,0,w\}}c[t',f_t\times \{d\},z].$$
Notice that the color assigned to $x_0$ can be $+\infty,0,1,\ldots,w$, but $f_t(x_0)<w$ means that $x_0$ can not be dominated $w$ times because $x_0$ will never appear in a bag for the rest of the algorithm. We know $c[t',f_t\times \{d\},z]\leq c[t',f_t\times \{d-1\},z]$ for integer $1\leq d<w$ and $z\in \{0,1,\ldots,L\}$. Thus we just consider the case the color assigned to $x_0$ is $\{+\infty,0,w\}$ while calculating $c[t,f_t,z]$. It takes $O((w+2)^{k+1}L)$ time for each forget node.

\textbf{Introduce node:} Suppose $t$ is an introduce node. Assume $B_{t'}=\{x_1,\ldots,x_{|B_t'|}\}$ and $B_t=B_{t'}\cup \{x_0\}$ for a $x_0\in V(G)$. For each  $f_t=(f_{t}(x_1),\ldots,f_{t}(x_{|B_{t'}|}),f_t(x_0))$, we have $f_t=f_{t'}\times \{f_t(x_0)\}$, where $f_{t'}=(f_{t}(x_1),\ldots,f_{t}(x_{|B_{t'}|}))$.
 Define $f'_{t'}=(f'_{t'}(x_1),\ldots,f'_{t'}(x_{|B_{t'}|}))$ such that for all $y\in B_{t'}$,

$$f'_{t'}(y)=\left\{
\begin{array}{ll}
+\infty  &   \mbox{if } f_{t'}(y)=+\infty ,\\
\max(0,f_{t'}(y)-1)  &  \mbox{if } f_{t'}(y)\neq +\infty \mbox{ and } y\in N(x_0) ,\\
f_{t'}(y)           &  \mbox{otherwise.}
\end{array}\right.$$
Then $f'_{t'}$ is a coloring of $B_{t'}$. The calculation of the evaluation index is shown as follows:

\n (1) $c[t,f_{t'}\times \{+\infty\},z]=\left\{
\begin{array}{ll}
-\infty  &   \mbox{if } |f^{-1}_{t'}(+\infty)|\geq z, \\
c[t',f'_{t'},z-1]+1+|N(x_0)\cap f_{t'}^{-1}(w)|           &  \mbox{otherwise;}
\end{array}\right.$


\n(2) $c[t,f_{t'}\times \{w\},z]=\left\{
\begin{array}{ll}
-\infty  &   \mbox{if } |N(x_0)\cap f^{-1}_{t'}(+\infty)|< w, \\
c[t',f_{t'},z]+1          &  \mbox{otherwise;}
\end{array}\right.$



\n (3) for $u=0,1,\ldots,w-1$,

$c[t,f_{t'}\times \{u\},z]=\left\{
\begin{array}{ll}
-\infty  &   \mbox{if } |N(x)\cap f^{-1}_{t'}(+\infty)|< u, \\
c[t',f_{t'},z]          &  \mbox{otherwise.}
\end{array}\right.$



If $f_t(x_0)=+\infty$, then $x_0$ is selected into $Dom(t,f_t,z)$. Since $f_t(y)=f'_{t'}(y)$ for $y\in B_{t'}$,  $Dom(t',f'_{t'},z-1)\cup\{x_0\}=Dom(t,f_t,z)$. Thus we have (1).

If $f_t(x_0)=u$ for  $0\leq u\leq w$, we need to check whether $|N(x_0)\cap f^{-1}_{t'}(+\infty)|\geq u$ (notice that there is no edge between $x_0$ and $y$ for all $y\in \cup_{t''\in T_{t'}}B_{t''}\setminus B_{t'}$). Thus (2) and (3) hold.

Given a coloring mapping $f_t$ and $z\in\{0,1,\ldots,L\}$, we need $O(k+1)$ time to calculate the number of neighbors selected into the $w$-dominating set. Then it takes $O((w+2)^{k+1}(k+L))$ time for each introduce node.

\textbf{Join node:} Suppose $t$ is a join node. Assume its child nodes are $t_1,t_2$ with $B_t=B_{t_1}=B_{t_2}$. For each coloring mapping $f_t$, we define $\{f_{t_1},f_{t_2}\}$ to be a good pair of $f_t$ as Theorem 3.1.


According to the definition of good pair, for a coloring mapping $f_t$, there exists a good pair $\{f_{t_1},f_{t_2}\}$ such that
$$Dom(t,f,z)=Dom(t_1,f_{t_1},z_1)\cup Dom(t_2,f_{t_2},z_2),$$
where $z=z_1+z_2-|f_t^{-1}(+\infty)|$ and $z_1,z_2\geq |f_t^{-1}(+\infty)|$. Since $(\cup_{t'\in T_{t_1}}B_{t'})\cap (\cup_{t'\in T_{t_2}}B_{t'})\subseteq B_t$, we have $Dom(t_1,f_{t_1},z_1)\cap Dom(t_2,f_{t_2},z_2)=f_{t_1}^{-1}(+\infty)=f_{t_2}^{-1}(+\infty)$. Then
$$|Dom(t,f,z)|=|Dom(t_1,f_{t_1},z_1)|+|Dom(t_2,f_{t_2},z_2)|-|f_{t}^{-1}(+\infty)|.$$
If $x\in dom(Dom(t,f_t,z)\rightarrow \cup_{t'\in T_{t}}B_{t'})\cup Dom(t,f_t,z)$, the following situations may occur: (i) $x\in B_t$ and $f_t(x)=+\infty$; (ii) $x\in B_t$ and $f_t(x)=w$; (iii) $x\notin B_t$ and $x\in Dom(t_i,f_{t_i},z_i)$ for an $i\in \{1,2\}$; (iv) $x\notin B_t$ and $|N(x)\cap Dom(t_i,f_{t_i},z_i)|\geq w$ for an $i\in \{1,2\}$. Then
$$\begin{array}{rcl}
& &|dom(Dom(t,f_t,z)\rightarrow\cup_{t'\in T_{t}}B_{t'})|+|Dom(t,f_t,z)|\\
&=&|dom(Dom(t_1,f_{t_1},z_1)\rightarrow\cup_{t'\in T_{t_1}}B_{t'})|+|Dom(t_1,f_{t_1},z_1)|\\
& &+|dom(Dom(t_2,f_{t_2},z_2)\rightarrow\cup_{t'\in T_{t_2}}B_{t'})|+|Dom(t_2,f_{t_2},z_2)|\\
& &-|f_{t}^{-1}(+\infty)|-|f^{-1}_{t_1}(w)|-|f^{-1}_{t_2}(w)|+|f^{-1}_{t}(w)|.
\end{array}$$
For disjoint sets $R_{+\infty},R_2,\ldots,R_{w}\subseteq B_t$, let $\overline{R}=R_{+\infty}\cup R_2\cup\ldots\cup R_w$ and $$\mathcal{F}_{\overline{R}}=\{f_t|f_t^{-1}(i)=R_i,i\in \{2,\ldots,w, +\infty\}\}.$$ 
For each $f_t\in \mathcal{F}_{\overline{R}}$, by the same argument as Theorem 3.1, we can think about that $f_t\in \mathcal{F}_{\overline{R}}$ is determined by $f_t^{-1}(1)$. Let $R_1\subset B_t\setminus \overline{R}$, we will rewrite $f_t\in \mathcal{F}_{\overline{R}}$ with $f_t^{-1}(1)=R_1$ by $f_t^{R_1}$. Given $R_1\subset B_t\setminus \overline{R}$ and $f_t^{R_1}\in \mathcal{F}_{\overline{R}}$. Let $f_t^{R_1}(x)=s_x$ for  $x\in B_t\setminus R_{+\infty}$. Denote
$$\mathcal{B}_{R_{1}}=\{\{f_{t_1},f_{t_2}\}|f^{-1}_{t_1}(+\infty)=f^{-1}_{t_2}(+\infty)=R_{+\infty},f_{t_1}(x)+f_{t_2}(x)-s'(x)=s_x, 2\le s_x\le w\}$$
and
$$\mathcal{C}_{z}=\{\{z_1,z_2\}|z_1+z_2-|f_t^{-1}(+\infty)|=z,z_1,z_2\geq |f_t^{-1}(+\infty)|\}.$$
Then $|\mathcal{B}_{R_{1}}|\le \prod\limits_{s=2}^{s=w}(s+1)^{|R_s|}$ and $|\mathcal{C}_z|\leq z+1$.
Then the evaluation index of $f_t$ is
$$c[t,f_t,z]=c[t,f_t^{R_1},z]=\max\limits_{\{f_{t_1}^{R_1},f_{t_2}^{R_1}\}\in \mathcal{B}_{R_{1}}}c[t,f_t^{R_1},f_{t_1}^{R_1},f_{t_2}^{R_1},z],$$
where
$$\begin{array}{rcl}
& &c[t,f_t^{R_1},f_{t_1}^{R_1},f_{t_2}^{R_1},z]\\
&=&\max\limits_{\substack{R^1_1\cup R^2_1=R_1\\R^1_1\cap R^2_1=\emptyset\\\{z_1,z_2\}\in\mathcal{C}_z}} \{c[t_1,f_{t_1},z_1]+c[t_2,f_{t_2},z_2]-|f_{t}^{-1}(+\infty)|-|f^{-1}_{t_1}(w)|-|f^{-1}_{t_2}(w)|+|f^{-1}_{t}(w)|\}.\end{array}$$

By Lemma 2.2, we can compute $c[t,f_t^{R_1},f_{t_1}^{R_1},f_{t_2}^{R_1}, z]$ for every $f_t\in \mathcal{F}_{\overline{R}}$, a good pair and a $z\in\{0,1,\ldots,L\}$ in $z2^{|B_t\setminus (R_{+\infty}\cup R_2\cup\ldots\cup R_w)|}|B_t|^{O(1)}$ time, since there are at most $z+1$ pairs $(z_1,z_2)\in \mathcal{C}_z$. By the similar discussion as Theorem 3.1, the total complexity of computing coloring mapping $f_t$ is

$$
\sum\limits_{z=0}^{L}z(\frac{(w+1)(w+2)}{2})^{|B_t|}|B_t|^{O(1)}=
(\frac{L(L+1)}{2})(\frac{(w+1)(w+2)}{2})^{|B_t|}|B_t|^{O(1)}.
$$
The total time spent for each join node is $(\frac{L(L+1)}{2})(\frac{(w+1)(w+2)}{2})^{k+1}(k+1)^{O(1)}$. Notice that the calculation of $c[t,f_t,z]$ of join node is more complicated than that of forget node and introduce node for integer $w,s\geq 1$, the total time for the algorithm is $(\frac{L(L+1)}{2})(\frac{(w+1)(w+2)}{2})^{k}k^{O(1)}n$.

\q

\vskip.2cm

If we set $w=1$, we get the following corollary.

{\bf Corollary 4.1 }
{\em Let $k,n$ be two positive integers and $G=(V,E)$ be a graph of order $n$ with $tw(G)=k$. The $L$-Max Dominating Set problem can be solved in time
$O(L^23^kk^{O(1)}n).$}

Corollary 4.1 improves the result proved by Roayaei et al. in \cite{RImpro} who gave a $4^{k}L^2n^{O(1)}$ time algorithm  for the $L$-Max Dominating Set problem.
\vskip.2cm

{\bf Acknowledgements~} This work is partially supported by the National Natural Science Foundation of China (Grant 11771247 \& 11971158) and  Tsinghua University Initiative Scientific Research Program.

\end{document}